\documentclass[11pt]{amsart}
\usepackage{amssymb,amsmath,amsthm,enumerate}

\newtheorem{thm}{Theorem}[section]
\newtheorem{lem}[thm]{Lemma}

\newtheorem{prop}[thm]{Proposition}
\newtheorem{conj}[thm]{Conjecture}

\theoremstyle{remark}
\newtheorem{exm}{Example}
\newtheorem{rem}[thm]{Remark}
\theoremstyle{definition}
\newtheorem{defi}[thm]{Definition}

\newcommand {\Zz} {\mathbb{Z}}
\newcommand {\Nz} {\mathbb{N}}
\newcommand {\Cz} {\mathbb{C}}
\newcommand {\Qz} {\mathbb{Q}}
\newcommand {\Rz} {\mathbb{R}}
\newcommand {\Sz} {\mathbb{S}}
\newcommand {\CF} {{\mathcal F}}

\newcommand{\BAR}{\overline}
\DeclareMathOperator{\SL}{SL}
\DeclareMathOperator{\GL}{GL}
\DeclareMathOperator{\Irr}{Irr}
\DeclareMathOperator{\Gal}{Gal}
\DeclareMathOperator{\diag}{diag}
\DeclareMathOperator{\w}{w}

\DeclareMathOperator{\PSL}{PSL}

\begin{document}

\title{Integral modular data and congruences}
\author{Michael Cuntz}
\address{Michael Cuntz, Universit\"at Kaiserslautern,
Postfach 3049, 67653 Kaiserslautern}
\email{cuntz@mathematik.uni-kl.de}

\begin{abstract}
We compute all fusion algebras with symmetric rational $S$-matrix
up to dimension 12. Only two of them may be used
as $S$-matrices in a modular datum:
the $S$-matrices of the quantum doubles of $\mathbb{Z}/2\mathbb{Z}$ and $S_3$.
Almost all of them satisfy a certain congruence which
has some interesting implications, for example for their degrees.
We also give explicitly an infinite sequence of modular data with rational
$S$- and $T$-matrices which are neither tensor products of smaller
modular data nor $S$-matrices of quantum doubles of finite groups.
For some sequences of finite groups (certain subdirect products of $S_3,D_4,Q_8,S_4$),
we prove the rationality of the $S$-matrices of their quantum doubles.
\end{abstract}


\maketitle

\section{Introduction}

The modular datum is a structure which appears
in physics as well as in representation theory.
It mainly consists of complex matrices $S$ and $T$ defining
a representation of the modular group $\SL_2(\Zz)$
where the matrix $S$ may be used to define structure
constants of a fusion algebra via the formula of Verlinde.
There is almost nothing known about how to classify
modular data, it is even unknown if there exist infinitely
many modular data in some dimension.

The fusion algebras associated to some special cases
have been investigated thoroughly, for example the affine
data coming from Kac-Moody algebras have completely been
classified by Gannon in \cite{tG2}. Another important
class of modular data is given by the quantum doubles
of finite groups (an overview may be found in \cite{CGR})
and there is also a version twisted by a $3$-cocycle
(see \cite{DVVV}).

A very interesting connection is the one to the representation theory
of finite groups of Lie type (see \cite{Lus1}).
Also, there is a new theory in which a new structure called
`Spetses' is associated to a complex reflection group
in a similar way as in Lusztig's theory (see \cite{gMSp}).
There exist modular data for the Spetses
(see for example \cite{MUG} for the imprimitive case).
Some of the corresponding fusion algebras have
negative structure constants. The analogous constructions
for finite Coxeter groups also yield Fourier matrices,
but they are not always symmetric (see \cite{mGgM}).
In fact, this shows that there
are many structures which are almost modular data in the
sense of Gannon \cite{tG} and therefore there is no
`correct' definition for a modular datum.

Here, we try to approach a classification of
a restricted case, namely the case in which the matrix
$S$ is rational. The number of such matrices $S$ in each dimension
is then finite, so we can compute all of them up to
a given dimension. It is very surprising that there
are only two rational $S$-matrices up to dimension $12$
(up to permutations of columns):
the $S$-matrices of the quantum doubles of $\Zz/2\Zz$ and $S_3$.

The symmetric rational matrices $S$ in dimension up to $12$ which
define fusion algebras but do not necessarily have an associated
$T$-matrix enjoy some very interesting properties. One
of them is a certain congruence property which we try
to understand in more details. We call such matrices
congruence multiplication tables.
We prove some relations between their degrees and
that they exist in all dimensions except $2,3,5$ and $6$.

We also give an infinite sequence of modular data with
rational $S$- and $T$-matrices. The smallest new case in
this sequence has dimension $94$ and is neither a tensor product
of smaller modular data nor an $S$-matrix of a quantum double
of a finite group (even its fusion algebra does not come
from a quantum double). This shows how far away we are from
a classification: there are still examples which
do not appear in the list of examples given by Gannon in \cite{tG}.

In the next section, we give some infinite sequences
of finite group modular data which are `integral' and not just
tensor products of smaller data.
The last section is a collection of open questions.

The appendix contains a complete list of all
symmetric rational Fourier matrices (see definition \ref{intmoddef})
up to dimension $12$ and a short explanation how it has been computed.

In dimension $12$ there is a congruence multiplication
table which has a $T$-matrix of infinite order. It
has negative structure constants.
In dimension $16$, there is also a congruence multiplication
table which is not the $S$-matrix of a quantum double and
has a $T$-matrix of infinite order. It has nonnegative
structure constants.

In dimension $18$ there exists a modular datum with
rational $S$-matrix.
It is the smallest datum known to the author whose
fusion algebra is not isomorphic to a fusion algebra
coming from a quantum double,
which is rational, has positive structure constants and
for which the group $\langle S,T\rangle$ is finite.
Up to $4$-th roots of unity, this is the datum belonging
to the $5$-th family (in the numbering of the package {\sc Chevie}, \cite{GHLMMP})
of the exceptional complex reflection group $G_8$
(the $8$-th group in the Shephard-Todd numbering (see \cite{gcSjaT})).

\section{Modular data}

Modular data are structures being studied in a large
variety of contexts. The most common definition is the following
(compare \cite{tG}):

\begin{defi}
Let $n\in\Nz$. By a {\bf modular datum} we mean a pair
$S,T\in\Cz^{n\times n}$ of matrices such that:
\begin{enumerate}
\item \label{modrep1} $S$ is unitary and symmetric: $S\bar S^T=1$, $S=S^T$,
\item $T$ is diagonal and of finite order,
\item $S_{i1}>0$ for $1\le i \le n$,
\item \label{modrep2} $S^2 = (ST)^3$,
\item \label{modrep3} $N_{ij}^l:=\sum_k \frac{S_{ki}S_{kj}\BAR{S_{kl}}}{S_{k1}} \in \Zz$
 for all $1\le i,j,l \le n$.
\end{enumerate}
\end{defi}

Because of (\ref{modrep1}), (\ref{modrep2}) and (\ref{modrep3}),
the matrices $S$ and $T$ define a representation of the modular group $\mbox{SL}_2(\Zz)$,
which is the reason why the datum above is called modular
(see for example \cite{tG}, 4 for a proof).

\begin{defi}
We will call a matrix $S$ satisfying the axioms (1), (3) and (5)
a {\bf Fourier matrix} as in \cite{gL}.
\end{defi}

The numbers $N_{ij}^l$ are called structure constants and define
an associative commutative $\Zz$-algebra which is a free $\Zz$-module
of finite rank (see for example \cite{mC}, 2.1.3).
Algebras coming from a Fourier matrix
or a character table via such structure constants
are usually called table algebras or fusion algebras
(complex conjugation on the set of columns corresponds to
an involution with certain properties).
Remark that we do not assume that the structure constants $N_{ij}^l$
are nonnegative, because we would exclude many interesting examples.
This is not so usual.

\begin{defi}
We call the algebra defined by the structure constants of
a Fourier matrix a {\bf fusion algebra}.
\end{defi}

Gannon (\cite{tG}) gives a list of open questions involving
modular data. The most important one is of course a classification.
However this is probably a very difficult task because there
are even much more fusion algebras than the long list of
examples given by Gannon in \cite{tG}.
It is also unknown, whether there exists a dimension $n$ for
which there are infinitely many modular data. Most likely
there are some infinite sequences that remain to be found.

\section{Integral modular data}

\subsection{Definition}

We would like to approach the first of these problems by restricting
to certain modular data.
The numbers $S_{ij}/S_{i1}$ of an arbitrary modular datum are eigenvalues
of matrices $N_i$ with integer entries (see Remark \ref{posbou}) and
are therefore algebraic integers.
So in particular if $S$ has
only rational entries, then the matrix $s=(s_{ij})_{i,j}$ defined by
\[ s_{ij} := \frac{S_{ij}}{S_{i1}} \]
lies in $\Zz^{n\times n}$. This leads us to the
following definition:

\begin{defi}\label{intmoddef} Let $n\in\Nz$.
An {\bf integral modular datum} consists of matrices
$s\in\Zz^{n\times n}$ and $T\in\Cz^{n\times n}$ such that:
\begin{enumerate}
\item $s_{i1} = 1$ for $1\le i \le n$, $\det(s)\ne 0$,
\item $s {s}^T=\diag(d_1,\ldots,d_n)=:d$, so $d_i=\sum_j s_{ji}^2$,
\item $\sqrt{d_j} s_{ij} = \sqrt{d_i} s_{ji}$ for $1\le i,j \le n$ \quad ($s$ is symmetrizable),
\item $N_{ij}^l:=\sum_k \frac{s_{ki} s_{kj} s_{kl}}{d_k} \in \Zz$ for $1\le i,j,l \le n$,
\item $T$ is diagonal,
\item $S^2 = (ST)^3$ where $S_{ij}:=s_{ij}/\sqrt{d_i}$.
\end{enumerate}
Here, $\sqrt{d_i}$ is chosen greater than $0$.
The matrix $s$ is the {\bf integral Fourier matrix} of the datum.
\end{defi}

The simplest example is the datum where $s$ and $T$ are both the
$1\times 1$ matrix $s=T=(1)$. This example is so trivial that it
will be ignored in the remainder of the paper.

It does not really make a difference if we use $s$ instead
of $S$. The matrix $s$ has the advantage, that it has only
integer entries, but it has the disadvantage that it is not
symmetric. Remark also that
we do not require that $T$ has finite order. This is because
there exist some interesting examples which satisfy all axioms
except the finiteness of the group $\langle S,T \rangle$.

So a modular datum $(S,T)$ with $S\in\Qz^{n\times n}$ is integral,
and of course, integral modular data with $T$ of finite order are modular data.
But there exist also integral modular data for which the
matrix $S$ has irrational entries:
\begin{exm}\label{moddatZ2}
Let $\zeta\in\Cz$ be a primitive $24$-th root of unity. Then
\[ S:=\frac{1}{\sqrt{2}}\begin{pmatrix} 1 & 1 \\ 1 & -1 \end{pmatrix}, \quad
   s:=\begin{pmatrix} 1 & 1 \\ 1 & -1 \end{pmatrix}, \quad
   T:=\begin{pmatrix} \zeta & 0 \\ 0 & -\zeta^7 \end{pmatrix}, \]
define an integral modular datum.
\end{exm}

\begin{rem}\label{posbou}
If the structure constants of $s$ are nonnegative, then
the theorem of Perron-Frobenius on nonnegative matrices
applies to the matrices $N_i:=(N_{ij}^k)_{j,k}$ for $1\le i\le n$:
if $s_k$ is the $k$-th row vector of $s$, then
\[ s_{ki} s_k^T = N_i s_k^T, \]
because the rows of $s$ are the one dimensional representations
of the fusion algebra.
So $s_k^T$ is an eigenvector of $N_i$ to the eigenvalue $s_{ki}$.
The theorem says that there exists an eigenvalue of maximal
absolute value and that it has an eigenvector with positive real
coordinates. Since $s_1$ is the only row with positive coordinates
(because the rows are orthogonal), this gives us the bound
\begin{equation}\label{posbouequ} |s_{ij}| \le s_{1j} \end{equation}
for all $i,j$.
Conversely, it does not follow from this bound that the structure
constants are nonnegative:
There exist Fourier matrices $s$ with
negative structure constants satisfying this bound.
\end{rem}

Some well known examples for integral modular data are given by
the quantum double of a finite group. The category of representations
of this algebra is a modular tensor category (see for example \cite{bBaK}, 3.2).
For example for the groups $C_2,S_3,S_4,D_4$ we get
indeed matrices $S$ with rational entries.

For given modular data $(S_1,T_1)$ and $(S_2,T_2)$, the tensor
product matrices $S_1\otimes S_2$, $T_1\otimes T_2$ form again
a modular datum. So from now on, we will only consider matrices
which do not decompose as tensor products. Note however that a decomposition
in a tensor product of matrices which may not be decomposed
further is not unique, see the example below.
\begin{defi}
We will call a modular datum {\bf indecomposable} if it may
not be written as the tensor product of two smaller
modular data.
\end{defi}

Also, we will consider modular data only up to a
{\it simultaneous} permutation of rows and columns fixing
the $1$.
Remark that for the fusion algebra, it makes no difference
if you permute rows or columns (independently).

\begin{defi}
Let $(s,T)$, $(s',T')$ be integral modular data of dimension $n$.
Then $(s,T)$ and $(s',T')$ are called {\bf isomorphic} if there
exists a permutation $\sigma\in S_n$ such that
\[ s_{ij}=s'_{\sigma(i)\sigma(j)}, \quad
   T_{ij}=T'_{\sigma(i)\sigma(j)} \quad \mbox{for all } 1\le i,j\le n.\]
\end{defi}

In the following lemmas, $(s,T)$ is an integral modular datum and
\[ s {s}^T=\diag(d_1,\ldots,d_n)=:d. \]

\begin{lem}\label{finnum}
For a given dimension $n$, there is only a finite number
of integral Fourier matrices.
\end{lem}
\begin{proof}
The symmetry of $S$ and \ref{intmoddef}, (1) give
\[ \frac{\sqrt{d_1}}{\sqrt{d_i}}=\frac{\sqrt{d_1}}{\sqrt{d_i}} s_{i1}=s_{1i}, \]
so $d_1 = \sum_i s_{1i}^2 = \sum_i \frac{d_1}{d_i}$, or
\begin{equation}\label{suminv}
\frac{1}{d_1}+\cdots+\frac{1}{d_n} = 1,
\end{equation}
which has a finite number of solutions $(d_1,\ldots,d_n)\in\Zz_{>0}^n$.
The $i$-th row vector of $s$ has integer entries and length $\sqrt{d_i}$.
So the number of integral Fourier matrices with given $d_1,\ldots,d_n$ is finite.
\end{proof}

In fact, this shows even more: since $\frac{\sqrt{d_1}}{\sqrt{d_i}} = s_{1i} \in \Zz$,
$d_i$ is a divisor of $d_1$ and $\frac{d_1}{d_i}\in\Zz$ is a square for all $1\le i\le n$.
Also, $\frac{d_i}{d_j}$ is a square in $\Qz$ for all $1\le i,j\le n$.

\begin{lem}
The integer $d_1$ is the least common multiple of $d_2,\ldots,d_n$.
\end{lem}

\begin{proof}
Because of equation (\ref{suminv}) of lemma \ref{finnum}, we may write
\[ \frac{1}{d_1}=1-\frac{a}{b}, \]
where $b$ is the least common multiple of $d_2,\ldots,d_n$.
But then $\frac{1}{d_1}=\frac{b-a}{b}$ and $b$ divides $d_1$
(because $\frac{d_1}{d_i}\in\Zz$). Thus $b=d_1$ and $b-a=1$.
\end{proof}

\begin{lem}\label{oddeven}
If $n$ is odd, then all $d_1,\ldots,d_n$ are squares.\\
If $n$ is even and $q$ is the square free part of $d_n$, then
$\frac{d_1}{q},\ldots,\frac{d_n}{q}$ are squares in $\Zz$.
\end{lem}

\begin{proof}
The determinant of $s$ is an integer and
$(\det s)^2 = d_1\cdots d_n$.
The preceding lemmas tell us that
$\frac{d_1}{d_2},\ldots,\frac{d_1}{d_n}$ are integers and squares.
So
\[ d_1^n = d_1 d_2\frac{d_1}{d_2} \cdots d_n\frac{d_1}{d_n} \]
is a square and therefore all $d_i$, $1\le i\le n$ are squares
if $n$ is odd.
If $n$ is even, this tells us that all square free parts of the $d_i$
are equal.
\end{proof}

Apart from the abelian groups, there is no character table of a finite
group which may be used as Fourier matrix in a modular datum
(because of the lack of symmetry). But in
analogy to character tables, $d_1$ is something like the order of a group
and $s_{11},\ldots,s_{1n}$ are the degrees of the irreducible characters.
\begin{defi}
We will call $d_1,\ldots,d_n$ the {\bf norms},
$s_{11},\ldots,s_{1n}$ the {\bf degrees} and $d_1$ the {\bf size}
of the modular datum.
\end{defi}

An integral Fourier matrix is not uniquely determined by
its degrees, for example:
\begin{exm}\label{example2}
\begin{tiny}
\[ \left( \begin{array}{ccccccccccccc}
1 & 1 & 1 & 1 & 2 & 2 & 2 & 4 & 4 & 4 & 8 & 8 & 8 \\
1 & 1 & 1 & 1 & 2 & 2 & 2 & 4 & 4 & 4 & 8 & -8 & -8 \\
1 & 1 & 1 & 1 & 2 & 2 & 2 & 4 & 4 & 4 & -8 & 8 & -8 \\
1 & 1 & 1 & 1 & 2 & 2 & 2 & 4 & 4 & 4 & -8 & -8 & 8 \\
1 & 1 & 1 & 1 & -2 & 2 & 6 & 0 & 0 & -4 & 0 & 0 & 0 \\
1 & 1 & 1 & 1 & 2 & 2 & 2 & -4 & -4 & 4 & 0 & 0 & 0 \\
1 & 1 & 1 & 1 & 6 & 2 & -2 & 0 & 0 & -4 & 0 & 0 & 0 \\
1 & 1 & 1 & 1 & 0 & -2 & 0 & 2 & -2 & 0 & 0 & 0 & 0 \\
1 & 1 & 1 & 1 & 0 & -2 & 0 & -2 & 2 & 0 & 0 & 0 & 0 \\
1 & 1 & 1 & 1 & -2 & 2 & -2 & 0 & 0 & 0 & 0 & 0 & 0 \\
1 & 1 & -1 & -1 & 0 & 0 & 0 & 0 & 0 & 0 & 0 & 0 & 0 \\
1 & -1 & 1 & -1 & 0 & 0 & 0 & 0 & 0 & 0 & 0 & 0 & 0 \\
1 & -1 & -1 & 1 & 0 & 0 & 0 & 0 & 0 & 0 & 0 & 0 & 0 \\
\end{array} \right) \]
\[ \left( \begin{array}{ccccccccccccc}
1 & 1 & 1 & 1 & 2 & 2 & 2 & 4 & 4 & 4 & 8 & 8 & 8 \\
1 & 1 & 1 & 1 & 2 & 2 & 2 & 4 & 4 & 4 & 8 & -8 & -8 \\
1 & 1 & 1 & 1 & 2 & 2 & 2 & 4 & 4 & 4 & -8 & 8 & -8 \\
1 & 1 & 1 & 1 & 2 & 2 & 2 & 4 & 4 & 4 & -8 & -8 & 8 \\
1 & 1 & 1 & 1 & 2 & 2 & 2 & 4 & -4 & -4 & 0 & 0 & 0 \\
1 & 1 & 1 & 1 & 2 & 2 & 2 & -4 & 4 & -4 & 0 & 0 & 0 \\
1 & 1 & 1 & 1 & 2 & 2 & 2 & -4 & -4 & 4 & 0 & 0 & 0 \\
1 & 1 & 1 & 1 & 2 & -2 & -2 & 0 & 0 & 0 & 0 & 0 & 0 \\
1 & 1 & 1 & 1 & -2 & 2 & -2 & 0 & 0 & 0 & 0 & 0 & 0 \\
1 & 1 & 1 & 1 & -2 & -2 & 2 & 0 & 0 & 0 & 0 & 0 & 0 \\
1 & 1 & -1 & -1 & 0 & 0 & 0 & 0 & 0 & 0 & 0 & 0 & 0 \\
1 & -1 & 1 & -1 & 0 & 0 & 0 & 0 & 0 & 0 & 0 & 0 & 0 \\
1 & -1 & -1 & 1 & 0 & 0 & 0 & 0 & 0 & 0 & 0 & 0 & 0 \\
\end{array} \right) \]
\end{tiny}
Remark that the first matrix defines negative structure
constants. The second one is a congruence multiplication table
(see \ref{conmultab}) with positive structure constants.
\end{exm}

It is a surprising fact that up to dimension $11$ there are no other indecomposable
integral modular data than example \ref{moddatZ2} and the Fourier matrix $s$
of the quantum double of $S_3$ (listed in the appendix).
To be precise, there are $3$ permutations of rows of $s$ which give
non isomorphic modular data, because they can not be obtained by a simultaneous
permutation of rows and columns.
In dimension $12$ there is (up to independent permutation of rows and columns)
exactly one integral Fourier matrix having a $T$-matrix; but it has negative structure constants
and the $T$-matrix has infinite order.

\subsection{Integral modular data with a unique norm}

Remember that $d_1$ is divisible by all $d_i$, $1\le i\le n$.

\begin{lem}\label{cong1}
We have
\[ \sum_i s_{ij}^3\delta_i \equiv
   \sum_i s_{ij}^2\delta_i \equiv
   \sum_i s_{ij}  \delta_i \equiv 0 \:\: (\mbox{\rm mod } d_1), \]
where $\delta_i:=d_1/d_i \in \Zz$.
\end{lem}
\begin{proof}
Apply \ref{intmoddef}, (1) to \ref{intmoddef}, (4).
\end{proof}

We would like to prove that an integral Fourier matrix $s$ with
a unique norm has only $\pm 1$ as entries. This seems to be true
in general:

\begin{conj}\label{equdiconj}
If $d_1=\ldots=d_n$, then $s\in\{\pm 1\}^{n\times n}$.
\end{conj}

Assuming the bound (\ref{posbouequ}) of remark \ref{posbou}, this
is easy to prove:
(remember that the bound holds when the structure constants are nonnegative).

\begin{lem}\label{equdi}
Assume that $|s_{ij}| \le s_{1j}$ for all $i,j$.
If $d_1=\ldots=d_n$, then $s\in\{\pm 1\}^{n\times n}$.
\end{lem}

\begin{proof}
If $d_1=\ldots=d_n$, then the matrix $s$ is symmetric and
$d_i=n$ for all $i$.
Then $s_{1i}=1$ for all $i$, so by our assumption, $|s_{ij}| \le 1$
for all $i,j$. But $d_i=\sum_{j=1}^n s_{ij}^2=n$, hence $s_{ij}=\pm 1$
for all $i,j$.
\end{proof}

\begin{prop}
The only dimensions $n$ for which there exist integral
modular data $(s,T)$ with $d_1=\ldots=d_n$ and $|s_{ij}| \le s_{1j}$ for all $i,j$
are powers of $2$. \\
The corresponding fusion algebras are isomorphic to
the group ring of
\[ (\Zz/2\Zz)^k, \quad  n=2^k. \]
\end{prop}
\begin{proof}
By lemma \ref{equdi}, integral modular data with
$d_1=\ldots=d_n=n$ and $|s_{ij}| \le s_{1j}$ for all $i,j$
have an $s$-matrix $s$ with entries
in $\{\pm 1\}$. The rows are orthogonal and the first row
is $s_1:=(1,\ldots,1)$.
Because of all our assumptions, the axiom for integrality of
the structure constants gives
\[ \frac{1}{n} \sum_k s_{ki} s_{kj} s_{kl} \in \Zz \]
for all $1\le i,j,l\le n$. If we denote the $i$-th column of $s$
by $s_i$ and componentwise multiplication of $s_i$ and $s_j$ by $s_i s_j$,
then this equation says
\[ s_i s_j \mbox{ is orthogonal to } s_l, \quad \mbox{or} \quad
s_i s_j = s_l \]
(Note that $s_i s_j=-s_l$ is impossible because $s_{1i}=s_{1j}=s_{1l}=1$).
For given $i,j$ there has to be an $l$ such that $s_i s_j = s_l$
because otherwise $s_i s_j$ would be orthogonal to all
elements of an orthogonal basis.\\
The first row is the unit element with respect to componentwise
multiplication, so the set of rows forms an abelian group and
all elements of this group have order $2$.
\end{proof}

Remark that there exist symmetric matrices with entries
$\pm 1$ and orthogonal rows for example in dimension $12$
and it is conjectured that they exist in all dimension $n=4k$,
$k\in\Nz$ (they are called symmetric Hadamard matrices).
But in general, the structure constants of the corresponding fusion algebras
are rational numbers.

This does not exclude the existence of non isomorphic
integral modular data with $d_1=\ldots=d_n$ in the same dimension;
for example:
\[
s_1 = \begin{pmatrix}
     1 & 1 & 1 & 1 \\
     1 & 1 &-1 &-1 \\
     1 &-1 & 1 &-1 \\
     1 &-1 &-1 & 1
\end{pmatrix}, \quad
s_2 = \begin{pmatrix}
     1 & 1 & 1 & 1 \\
     1 & 1 &-1 &-1 \\
     1 &-1 &-1 & 1 \\
     1 &-1 & 1 &-1
\end{pmatrix}.
\]
Each of $s_1, s_2$ have 12 different solutions for the matrix $T$.

These matrices give also an example for non uniqueness of
a complete decomposition in tensor products:\\
Let $s$ be the matrix of example \ref{moddatZ2} (the character table
of $\Zz/2\Zz$), so $s_2 \cong s\otimes s$. Then
\[ s\otimes s_1 \cong s\otimes s_2 \]
although $s_1\ncong s_2$ and $s_1$ is not a tensor product
of smaller symmetric matrices of dimension greater than $1$.

\subsection{Zeros in Fourier matrices}

\begin{lem}
If an integral modular datum $(s,T)$ has two different
norms, then there is a zero in the matrix $s$.
\end{lem}
\begin{proof}
If there are two different norms, then there is
an entry of $s$ unequal to $\pm 1$, say $s_{kl}$.
Suppose that all entries are non zero. Then
all norms are greater or equal to $n$ and
because $s_{kl}\ne\pm 1$, $d_k>n$. Then
\[ \frac{1}{d_1}+\cdots+\frac{1}{d_n} < 1, \]
in contradiction to equation (\ref{suminv}) of lemma \ref{finnum}.
\end{proof}

Looking at the matrices of the appendix, one might be led to
believe that an integral Fourier matrix has a $T$-matrix only
if there is no zero on the diagonal. This is wrong:
for example the $S$-matrix of the quantum double of the wreath
product $Z_2\wr S_4$ (the group of monomial $4\times 4$-matrices
with entries $\pm 1$) is integral and has a zero on the diagonal.
It is a $325\times 325$-matrix.

\section{Congruence multiplication tables}\label{conmultab}

\subsection{Congruences modulo $\sqrt{d_1}$}

The computation of all integral Fourier matrices up to
dimension $12$ provides some interesting observations.
One of the most remarkable is that most of them satisfy
a certain congruence:

\begin{defi}
Let $s$ be an integral Fourier matrix, i.e.
$s\in\Zz^{n\times n}$ and
\begin{enumerate}
\item $s_{i1} = 1$ for $1\le i \le n$, $\det(s)\ne 0$,
\item $s {s}^T=\diag(d_1,\ldots,d_n)$,
\item $\sqrt{d_j} s_{ij} = \sqrt{d_i} s_{ji}$ for $1\le i,j \le n$,
\item $N_{ij}^l:=\sum_k \frac{s_{ki} s_{kj} s_{kl}}{d_k} \in \Zz$ for $1\le i,j,l \le n$.
\end{enumerate}
If in addition $s$ satisfies
\begin{enumerate}
\item[(5)] $d_1$ is a square, $w:=\sqrt{d_1}\in \Nz$,
\item[(6)] $g_i s_{ij} \equiv g_i g_j \:\: (\mbox{\rm mod } w)$ for all $1\le i,j \le n$,
\end{enumerate}
where $g_i:=s_{1i}$ are the degrees of $s$, then we say that
$s$ is a {\bf congruence multiplication table}.
Sometimes, we will also need the bound (see \ref{posbou})
\begin{enumerate}
\item[(7)] $|s_{ij}| \le g_i$ for all $1\le i,j \le n$.
\end{enumerate}
\end{defi}

Define $\Sz_{ij}:=s_{ij} s_{1i}$ for $1\le i,j \le n$, so $\Sz = wS$.
Then the congruence (6) takes the following form:
\[ \Sz_{ij} \equiv g_i g_j \:\: (\mbox{\rm mod } w). \]
Because $d_1$ is the largest norm, the entries of $\Sz$ are bounded
by $w$,
\[ -w< \Sz_{ij} <w.\]
So in a congruence multiplication table with given degrees $g_1,\ldots,g_n$,
there are at most two possible values for each entry.
For a congruence multiplication table $s$, we will always denote
the above symmetric matrix by $\Sz$.

\begin{exm}
Since the degrees divide $w$, one could ask if $g_i g_j$
is always divisible by $\Sz_{ij}$. This is true only up to
dimension $12$:
\[ s = \begin{tiny}\left( \begin{array}{ccccccccccccc}
1 & 1 & 1 & 1 & 2 & 8 & 8 & 8 & 10 & 10 & 20 & 20 & 20 \\
1 & 1 & 1 & 1 & 2 & 8 & 8 & 8 & 10 & 10 & 20 & -20 & -20 \\
1 & 1 & 1 & 1 & 2 & 8 & 8 & 8 & 10 & 10 & -20 & 20 & -20 \\
1 & 1 & 1 & 1 & 2 & 8 & 8 & 8 & 10 & 10 & -20 & -20 & 20 \\
1 & 1 & 1 & 1 & 2 & 8 & 8 & 8 & -10 & -10 & 0 & 0 & 0 \\
1 & 1 & 1 & 1 & 2 & 3 & -2 & -2 & 0 & 0 & 0 & 0 & 0 \\
1 & 1 & 1 & 1 & 2 & -2 & 3 & -2 & 0 & 0 & 0 & 0 & 0 \\
1 & 1 & 1 & 1 & 2 & -2 & -2 & 3 & 0 & 0 & 0 & 0 & 0 \\
1 & 1 & 1 & 1 & -2 & 0 & 0 & 0 & -2 & 2 & 0 & 0 & 0 \\
1 & 1 & 1 & 1 & -2 & 0 & 0 & 0 & 2 & -2 & 0 & 0 & 0 \\
1 & 1 & -1 & -1 & 0 & 0 & 0 & 0 & 0 & 0 & 0 & 0 & 0 \\
1 & -1 & 1 & -1 & 0 & 0 & 0 & 0 & 0 & 0 & 0 & 0 & 0 \\
1 & -1 & -1 & 1 & 0 & 0 & 0 & 0 & 0 & 0 & 0 & 0 & 0 \\
\end{array} \right) \end{tiny}\]
is a congruence multiplication table.
\end{exm}

Remember that $(7)$ follows if the structure constants are nonnegative.
Up to dimension 12, there are 28 integral Fourier matrices with nonnegative
structure constants and only two of them are not
congruence multiplication tables, the character tables
of the groups $\Zz/2\Zz$ and $(\Zz/2\Zz)^3$. For these two
matrices, the normalized $S$-matrix is not rational, because
the norms are not squares.

In particular, the Fourier matrix of the quantum double of
$S_3$ is a congruence multiplication table.
It seems as if dimension $12$ is too small to find more
complicated examples, because experiments show that bigger
indecomposable Fourier matrices coming from quantum doubles
do not satisfy this congruence property anymore.

Also, no character table of a group $(\Zz/2\Zz)^k$, $k>3$ is a congruence
multiplication table. In this case we have $w=2^{\frac{k}{2}}>2$
and the entries of the matrix $\Sz$ are $\pm 1$. But then
$-1 \not\equiv 1 \:\: (\mbox{\rm mod } w)$, so (6) is not
satisfied.

However, we may also define $w$ to be the largest integer
such that $(6)$ is satisfied. For example the Fourier
matrix of the quantum double of $S_4$ has the congruence
property for $w=\frac{\sqrt{d_1}}{2}$.

But let us first consider the case $w=\sqrt{d_1}$.

\begin{prop}\label{congrest}
Let $\Sz$ be a congruence multiplication table
and assume also the bound $(7)$. Then
\begin{eqnarray*}
g_i g_j < \frac{w}{2} & \Rightarrow & \Sz_{ij}=g_i g_j, \\
w | g_i g_j & \Rightarrow & \Sz_{ij}=0, \\
g_i g_j < w & \Rightarrow & \Sz_{ij}\ge -\frac{w}{2},
\end{eqnarray*}
for all $1\le i,j \le n$.
\end{prop}

\begin{proof}
Assume that $g_i g_j < \frac{w}{2}$. Then the two integers
congruent to $g_i g_j$ modulo $w$ in the interval $(-w,w)$ are
$g_i g_j$ and $g_i g_j-w$. The case $\Sz_{ij}=g_i g_j-w$ would
imply
\[ g_i g_j = w-|\Sz_{ij}| \ge w-g_i g_j, \]
because $\Sz_{ij}<0$ and so $2g_i g_j \ge w$ which we excluded
at the beginning.
The second and third implications are obvious.
\end{proof}

The first two implications of the last proposition
say that a big part of the matrix is given by the degrees.
An interesting problem is to characterize the degrees
appearing in congruence multiplication tables. For example,
experiments show that there are always at least 2 degrees
equal to 1. We have no proof for this, but we know:

\begin{prop}
Let $s$ be a congruence multiplication table.
Then the number $m$ of degrees of $s$ equal to $1$ is less or
equal to $4$.

If $m>1$, then there is at least one degree equal to
$\frac{w}{2} = \sqrt{\frac{d_1}{4}}$
(or norm equal to $4$).

If in addition the structure constants are nonnegative,
then $m\in\{ 1,2,4 \}$.
\end{prop}
\begin{proof}
We may assume without loss of generality
that $g_1\le \ldots \le g_n$.
If there are two degrees equal to $1$ then $g_1=g_2=1$
and the first two rows of $s$ will have norm $d_1$.
The largest possible degree is $\frac{w}{2}$,
because it should divide $w$ and the
sum of the squares of the degrees is $w^2=d_1$.
Let $r\in\Zz$ be maximal such that $g_{n-r+1},\ldots,g_n=\frac{w}{2}$,
so $r$ is the number of degrees equal to $\frac{w}{2}$.
By proposition \ref{congrest}, the first $n-r$ entries of $s$ in
the second row are determined by the degrees:
\[ s_{2i} = g_i \quad \mbox{for } i\le n-r. \]
Now the second row has to be orthogonal to the first one,
so $r>0$ and $d_n = 4$. But the sum of the inverse of
the norms must be one, so $r\le 4$ and if $n\ne 4$ then
$r<4$.

Since by $(6)$ the values at the end of the second row
are congruent $\frac{w}{2}$ modulo $w$, the number
of possibilities for the second row to be orthogonal
to the first is at most $3$ (if $r=3$).
We get, that the number $m$ of degrees of $s$ equal to $1$
is at most $4$.

If the structure constants are nonnegative, then the
set of columns with degree $1$ forms a group with respect
to componentwise multiplication, so the case $m=3$ is excluded.
\end{proof}

\begin{prop}
Let $s$ be a congruence multiplication table of dimension $n$
and with degrees $g_1,\ldots,g_n$. Then
\[ \sum_{\stackrel{j=1}{2 g_i g_j< w}}^n g_j^2 \le \sum_{\stackrel{j=1}{2 g_i g_j\ge w}}^n g_j^2 \]
for all $1\le i\le n$.
\end{prop}
\begin{proof}
Use the first implication of proposition \ref{congrest} together
with the orthogonality between the first and the $i$-th row.
\end{proof}

\begin{prop}\label{oldcmtp6}
If $s$ is an integral Fourier matrix of dimension $n$
and of size $d_1$ which is a square,
then
\[ \tilde s:=\begin{pmatrix}
1 & 1 & 1 & 1 & 2s_{12} & \cdots & 2s_{1n} & b & b & b \\
1 & 1 & 1 & 1 & 2s_{12} & \cdots & 2s_{1n} & b & -b & -b \\
1 & 1 & 1 & 1 & 2s_{12} & \cdots & 2s_{1n} & -b & b & -b \\
1 & 1 & 1 & 1 & 2s_{12} & \cdots & 2s_{1n} & -b & -b & b \\
1 & 1 & 1 & 1 & 2s_{22} & \cdots & 2s_{2n} & 0 & 0 & 0 \\
\vdots &   &   & \vdots & \vdots &  & \vdots & \vdots &  & \vdots \\
1 & 1 & 1 & 1 & 2s_{n2} & \cdots & 2s_{nn} & 0 & 0 & 0 \\
1 & 1 & -1 & -1 & 0 & \cdots & 0 & 0 & 0 & 0 \\
1 & -1 & 1 & -1 & 0 & \cdots & 0 & 0 & 0 & 0 \\
1 & -1 & -1 & 1 & 0 & \cdots & 0 & 0 & 0 & 0 \\
\end{pmatrix} \]
where $b:=2\sqrt{d_1}$
is an integral Fourier matrix of dimension
$n+6$.

If $s$ is a congruence multiplication table, then
$\tilde s$ is one as well.
\end{prop}
\begin{proof}
Orthogonality of the rows of $\tilde s$ is an easy verification.
The size of $\tilde s$ is
\[ \tilde d_1 = 4+4(d_1-1)+12 d_1 = 16 d_1 \]
so the symmetry is obvious as well.
For the integrality of the structure constants, we consider all
componentwise products of columns; they evidently always yield linear
combinations of the columns with integer coefficients.
The second assertion is trivial.
\end{proof}

\begin{rem}
The contruction of proposition \ref{oldcmtp6} is something
like an extension of an integral Fourier matrix
by $s_2\otimes s_2$ where $s_2$ is the
character table of the cyclic group $\Zz/2\Zz$.
It is possible to generalize this to an extension
by $s_e\otimes s_e$ where $s_e$ is the
character table of the cyclic group $\Zz/e\Zz$, $e\in\Nz$.
The resulting matrix has $2(e^2-1)$ more columns and
is of course not rational anymore (if $e>2$). But the structure
constants remain rational integers.
\end{rem}

For a given matrix $s$, we will denote the matrix $\tilde s$ of proposition
\ref{oldcmtp6} obtained from $s$ by $s^{+6}$.

\begin{exm}
The matrix
\begin{tiny}
\[ s:=\left( \begin{array}{ccccccccccc}
1 & 1 & 1 & 3 & 3 & 3 & 3 & 5 & 5 & 6 & 10 \\
1 & 1 & 1 & 3 & 3 & 3 & 3 & 5 & -10 & 6 & -5 \\
1 & 1 & 1 & 3 & 3 & 3 & 3 & -10 & 5 & 6 & -5 \\
1 & 1 & 1 & -2 & 3 & -2 & -2 & 0 & 0 & 1 & 0 \\
1 & 1 & 1 & 3 & -2 & -2 & -2 & 0 & 0 & 1 & 0 \\
1 & 1 & 1 & -2 & -2 & -2 & 3 & 0 & 0 & 1 & 0 \\
1 & 1 & 1 & -2 & -2 & 3 & -2 & 0 & 0 & 1 & 0 \\
1 & 1 & -2 & 0 & 0 & 0 & 0 & -1 & -1 & 0 & 1 \\
1 & -2 & 1 & 0 & 0 & 0 & 0 & -1 & -1 & 0 & 1 \\
1 & 1 & 1 & \frac{1}{2} & \frac{1}{2} & \frac{1}{2} & \frac{1}{2} & 0 & 0 & -\frac{3}{2} & 0 \\
1 & -\frac{1}{2} & -\frac{1}{2} & 0 & 0 & 0 & 0 & \frac{1}{2} & \frac{1}{2} & 0 & -\frac{1}{2} \\
\end{array} \right) \]
\end{tiny}
is not an integral Fourier matrix, because its structure constants are
not all integers (and because its entries are not all integers).
But the construction of proposition \ref{oldcmtp6}
applied to $s$ yields a $17\times 17$-matrix which is a
congruence multiplication table.
\end{exm}

\begin{thm}
There exist congruence multiplication tables
in all dimensions except $2,3,5,6$.
\end{thm}
\begin{proof}
In the appendix, we give congruence multiplication
tables for the dimensions $n=4,7,8,9,10,11$ and $12$. By proposition
\ref{oldcmtp6}, they exist in all dimensions greater than $12$.
To exclude the cases $2,3,5$ and $6$ we have computed all integral
Fourier matrices in these dimensions.
\end{proof}

\subsection{The modulus of an integral Fourier matrix}

Now we consider the case where $w$ is not necessarily equal to $\sqrt{d_1}$.

\begin{defi}
Let $s$ be an integral Fourier matrix.
Then we call the greatest integer $w\in\Zz_{>0}$ such that
\[ g_i s_{ij} \equiv g_i g_j \:\: (\mbox{\rm mod } w) \mbox{ for all } 1\le i,j \le n, \]
where $g_i:=s_{1i}$ are the degrees of $s$,
the {\bf modulus} $\w(s):=w$ of $s$.
We set $\w((1)):=\infty$,
$a\equiv b \:\: (\mbox{\rm mod } \infty) :\Leftrightarrow a=b$,
and $\gcd(a,\infty):=a$ for all $a,b \in \Zz$.
\end{defi}

\begin{prop}
If $s$ and $s'$ are integral Fourier matrices, then
\[ \w(s\otimes s') = \gcd(\w(s),\w(s')). \]
\end{prop}

\begin{proof}
Let $n,n'$ be the dimensions of $s,s'$.
The rows and columns of the tensor product $s\otimes s'$ are indexed by pairs
$(i_1,i_2)$, $1\le i_1 \le n$, $1\le i_2 \le n'$ and an entry is
\[ (s\otimes s')_{(i_1,i_2),(j_1,j_2)} = s_{i_1 j_1}s'_{i_2 j_2}. \]
Then
\[ g_{i_1} g'_{i_2} s_{i_1 j_1}s'_{i_2 j_2} \equiv g_{i_1} g'_{i_2}g_{j_1} g'_{j_2}
\:\: (\mbox{\rm mod } \w(s\otimes s')) \]
and this holds in particular for $i_1=j_1=1$ and $i_2=j_2=1$
from which follows that $\w(s\otimes s')$ divides $\tilde w:=\gcd(\w(s),\w(s'))$.
On the other hand,
\[ g_i s_{ij} \equiv g_i g_j \:\: (\mbox{\rm mod } \tilde w) \mbox{ for all } 1\le i,j \le n, \]
\[ g'_i s'_{ij} \equiv g'_i g'_j \:\: (\mbox{\rm mod } \tilde w) \mbox{ for all } 1\le i,j \le n', \]
so $g_{i_1} g'_{i_2} s_{i_1 j_1}s'_{i_2 j_2} \equiv g_{i_1} g'_{i_2}g_{j_1} g'_{j_2}
\:\: (\mbox{\rm mod } \tilde w)$ and hence $\tilde w \le \w(s\otimes s')$.
\end{proof}

The set $\CF$ of integral Fourier matrices up to isomorphisms
forms a monoid together with tensor product. The preceding proposition
states that
\[ \w : (\CF,\otimes) \rightarrow (\Nz\cup\{\infty\},\gcd), \quad s \mapsto \w(s)  \]
is a monoid homomorphism.

\section{Integral modular data with two different degrees}\label{infint}

\subsection{Some general observations}
In this section we will try to find out how integral Fourier
matrices with two different degrees look like.
If the structure constants are positive, then the columns of
degree $1$ form a group isomorphic to a direct product
of cyclic groups $\Zz/2\Zz$, so the number of degrees equal
to $1$ is a power of $2$, say $2^m$, $m\in\Zz_{\ge 0}$.
If there are only two different values $1$ and $g$ for the degrees,
then $d_1 = 2^m+q g^2$ for some $q\in\Zz$. This implies
that $g^2$ divides $2^m$, because $g^2$ divides $d_1$.
So $g=2^k$ for some $k\in\Zz$ with $m\ge 2k$.

Let $c$ be the character table of $(\Zz/2\Zz)^m$.
Then the first $2^m$ columns of the matrix $s$ build
a $(2^m+q)\times 2^m$-matrix.
Let $a_1,\ldots,a_{2^m}$ be the multiplicities of the
rows of this matrix; so the $i$-th row of $c$ appears
$a_i$ times.

Let us first assume that $a_1=2^m$.
This holds for example in the infinite sequence of modular
data of section \ref{infintser}.
Then the symmetry and orthogonality of the Fourier matrix give
\[ a_1+\sum_{i=2}^{2^m} c_{ij} a_i g^2 = 0 \]
for all $1<j\le 2^m$ and
\[ a_1+\sum_{i=2}^{2^m} a_i g^2 = 2^m+q 2^{2k} \]
or equivalently
\[ (2^{m-2k},a_2,\ldots,a_{2^m}) c = (2^{m-2k}+q,0,\ldots,0). \]
But $c^{-1}=2^{-m}c$ and the first row of $c$ has only
$1$'s, so
\[ q=(2^m-1) 2^{m-2k}, \quad a_2 = \ldots = a_{2^m} = 2^{m-2k}. \]

Now we return to the general case where $a_1$ is arbitrary.
Let $0 \le \nu < 2^m$ be such that
\[ \sum_{i=1}^\nu a_i \le 2^m, \quad -2^m+\sum_{i=1}^{\nu+1} a_i =: \delta>0. \]
Then
\[ (\sum_{i=1}^\nu a_i c_{ij}) + (a_{\nu+1}-\delta+\delta g^2) c_{\nu+1,j}
+ \sum_{i=\nu+2}^{2^m} a_i c_{ij} g^2 = \delta_{ij} (2^m+qg^2). \]
By the same argument as above, we get
\begin{eqnarray*}
a_1,\ldots,a_\nu & = & (2^m+qg^2) 2^{-m}, \\
a_{\nu+1} & = & (2^m+qg^2) 2^{-m}+\delta-\delta g^2, \\
a_{\nu+2},\ldots,a_{2^m} & = & (2^{m-2k}+q)2^{-m}.
\end{eqnarray*}
Because $a_1,\ldots,a_\nu$ are integers, $2^{m-2k}$ divides $q$
say $q=q' 2^{m-2k}$. By the third equality, $2^{2k}$ divides $1+q'$,
so $2^{2k}\le 1+q'$ and therefore
\[ \nu \le \frac{2^m}{a_\nu} = \frac{2^{2m}}{2^m+qg^2}
= \frac{2^{m}}{1+q'} \le 2^{m-2k}. \]
Unfortunately, in a certain sense, there is no better bound:
\begin{exm}
The matrix
\begin{tiny}
\[ \left( \begin{array}{cccccccccccccc}
1 & 1 & 1 & 1 & 1 & 1 & 1 & 1 & 2 & 2 & 2 & 2 & 2 & 2 \\
1 & 1 & 1 & 1 & 1 & 1 & 1 & 1 & 2 & 2 & -2 & -2 & -2 & -2 \\
1 & 1 & 1 & 1 & 1 & 1 & 1 & 1 & -2 & -2 & 2 & 2 & -2 & -2 \\
1 & 1 & 1 & 1 & 1 & 1 & 1 & 1 & -2 & -2 & -2 & -2 & 2 & 2 \\
1 & 1 & 1 & 1 & -1 & -1 & -1 & -1 & 2 & -2 & 2 & -2 & 2 & -2 \\
1 & 1 & 1 & 1 & -1 & -1 & -1 & -1 & 2 & -2 & -2 & 2 & -2 & 2 \\
1 & 1 & 1 & 1 & -1 & -1 & -1 & -1 & -2 & 2 & 2 & -2 & -2 & 2 \\
1 & 1 & 1 & 1 & -1 & -1 & -1 & -1 & -2 & 2 & -2 & 2 & 2 & -2 \\
1 & 1 & -1 & -1 & 1 & 1 & -1 & -1 & 0 & 0 & 0 & 0 & 0 & 0 \\
1 & 1 & -1 & -1 & -1 & -1 & 1 & 1 & 0 & 0 & 0 & 0 & 0 & 0 \\
1 & -1 & 1 & -1 & 1 & -1 & 1 & -1 & 0 & 0 & 0 & 0 & 0 & 0 \\
1 & -1 & 1 & -1 & -1 & 1 & -1 & 1 & 0 & 0 & 0 & 0 & 0 & 0 \\
1 & -1 & -1 & 1 & 1 & -1 & -1 & 1 & 0 & 0 & 0 & 0 & 0 & 0 \\
1 & -1 & -1 & 1 & -1 & 1 & 1 & -1 & 0 & 0 & 0 & 0 & 0 & 0 \\
\end{array} \right) \]
\end{tiny}
is the tensor product of the matrix of dimension $2$ and the
matrix of dimension $7$ of the appendix. We have
\[ m=3, \quad k=1, \quad q=5, \quad \nu=2, \quad \delta=0,
\quad a_1=a_2=4, \quad a_3,\ldots,a_8=1 .\]
Remark also that this is not a congruence multiplication table.
\end{exm}

\subsection{An infinite sequence of integral modular data
with rational $T$-matrix and different degrees}\label{infintser}

The modular datum given by the quantum double of
the dihedral group with $8$ elements has the above properties.
By permuting two columns in the Fourier matrix we get
a Fourier matrix which has a rational $T$-matrix
of finite order (this is excluded in $T$-matrices
of quantum doubles of non abelian groups, because the order of
$T$ is exactly the exponent of the corresponding group).
We generalize this and get the following sequence.

\begin{prop}\label{infintexpl}
Let $k\in\Zz$,
\[ m:=2k+1, \quad g:=2^k, \quad q:=2(2^m-1), \quad a:=2^m, \]
and
\[ c:=\begin{pmatrix} 1 & 1 \\ 1 & -1 \end{pmatrix}
\otimes \cdots \otimes \begin{pmatrix} 1 & 1 \\ 1 & -1 \end{pmatrix} \]
be the $m$-th tensor power of the character table of $\Zz/2\Zz$.
Then
\[ S:=\frac{1}{a} \begin{pmatrix}
1 & \cdots & 1 & gc_{21} & gc_{21} & \cdots & gc_{a1} & gc_{a1} \\
\vdots & & \vdots & \vdots & \vdots & & \vdots & \vdots \\
1 & \cdots & 1 & gc_{2a} & gc_{2a} & \cdots & gc_{aa} & gc_{aa} \\
gc_{21} & \cdots & gc_{2a} & g^2 & -g^2 & & & \\
gc_{21} & \cdots & gc_{2a} & -g^2 & g^2 & & 0 & \\
\vdots & & \vdots & & & \ddots & & \\
gc_{a1} & \cdots & gc_{aa} & & 0 & & g^2 & -g^2  \\
gc_{a1} & \cdots & gc_{aa} & & & & -g^2 & g^2  \\
\end{pmatrix} \]
and
\[ T=\diag(\underbrace{1,\ldots,1}_{a},\underbrace{1,-1,\ldots,1,-1}_{2(a-1)}) \]
define a modular datum of dimension $2(3\cdot 2^{2k}-1)$
with rational matrices $S$ and $T$.

The structure constants are nonnegative and their greatest value is $2^{k-1}$.
\end{prop}

\begin{proof}
Orthogonality of the rows of $S$ is an easy verification
(use several times the fact that $c$ is orthogonal).

We want to prove that the structure constants are nonnegative rational integers.
For this, consider the matrix $s$ which is obtained by dividing each
row of $S$ by its first entry. Then the fusion algebra is the $\Zz$-lattice
spanned by the columns $v_1,\ldots,v_a,w_1,\ldots,w_{2a-2}$ of $s$ and
multiplication is componentwise. Products of the form $v_i v_j$ or
$v_i w_j$ give columns of $s$, so the only interesting cases are
$w_i w_j$. Without loss of generality, we may restrict to the three cases
$w_1 w_1$, $w_1 w_2$ and $w_1 w_3$. One easily sees that
\[ w_1 \cdot w_1 = \sum_{\substack{i=1 \\ s_{a+1,i}=1}}^a v_i, \quad
w_1 \cdot w_2 = \sum_{\substack{i=1 \\ s_{a+1,i}=-1}}^a v_i, \]
and that
\[ w_1 \cdot w_3 = 2^{k-1}(w_j+w_{j+1}) \quad \mbox{ for some } 1\le j \le 2a-2. \]
It remains to verify $(ST)^3 = 1$. A tiring calculation gives
$(ST)^2=TS$, which suffices because $T^2=S^2=1$.
\end{proof}

\begin{rem}
The subgroup $\langle S,T \rangle\le \GL_{2(3\cdot 2^{2k}-1)}(\Qz)$ spanned by the matrices
$S$ and $T$ is isomorphic to the symmetric group $S_3$.
\end{rem}

The proof of the following proposition for arbitrary $k$ has been suggested
by the referee.

\begin{prop}\label{referee1}
The fusion algebras of the $S$-matrices of Proposition \ref{infintexpl} for $k>1$
are not isomorphic to the fusion algebras associated to quantum doubles
of finite groups.
\end{prop}
\begin{proof}
Equation \ref{quantS} of section \ref{infintm} gives a definition
for the matrices $S$ and $T$. They are indexed by pairs $(\bar g,\chi)$,
where $\bar g$ is a conjugacy class of $G$ and $\chi$ is an irreducible
character of the centralizer of $g$ in $G$. In particular, the column
of $S$ indexed by $(\bar e,1)$ is the Perron-Frobenius vector,
this is the unique column of $S$ with positive real entries.
These entries are
\[ S_{(\bar h,\chi),(\bar e,1)} = \frac{\deg(\chi)}{|C_G(h)|} \]
and if our $S$-matrix was associated to the quantum double of a finite group $G$,
then these would be the values
$\frac{1}{a},\ldots,\frac{1}{a},\frac{g}{a},\ldots,\frac{g}{a}$,
where $\frac{1}{a}$ occurs $a$ times and $\frac{g}{a}$ occurs $2(a-1)$ times.
It follows immediately that $G$ has $a$ elements and is nonabelian.

Now if $h$ does not lie in the center $Z$ of $G$, then its centralizer
will have less than $a$ elements and since $S_{(\bar h,1),(\bar e,1)}=|C_G(h)|^{-1}$,
we get $|C_G(h)|=\frac{a}{g}$. This also implies that the centralizer is abelian
because all degrees of irreducible characters will be $1$.
Hence the $a$ values $\frac{1}{a}$ correspond to pairs $(\bar h,\chi)$ with $h\in Z$
and $\deg(\chi)=1$. So the number of $1$-dimensional irreducible representations of
$G$ is $a/|Z|$. Further, all irreducible representations of $G$ have dimension
$1$ or $g$.

Denote by $n_g$ the number of irreducible representations of $G$ of dimension
$g$. Then $a=a/|Z|+g^2 n_g$, i.e. $|Z|=2$ and $n_g=1$. This implies that
the number of conjugacy classes of $G$ is $2^{2k}+1$ and therefore we have
$2^{2k}-1$ classes which are not in $Z$. Finally we obtain
\[ a=2^{2k+1}=|G|=|Z|+(2^{2k}-1)\frac{a}{2^{k+1}}=2+2^{3k}-2^{k} \]
which is impossible for $k>1$.
\end{proof}

\section{Some infinite sequences of integral finite group modular data}\label{infintm}

For a finite group $G$, the $S$- and $T$-matrices of the quantum
double of $G$ are indexed by pairs $(\bar g,\chi)$,
where $\bar g$ is a conjugacy class of $G$ and
$\chi$ is an irreducible character of the centralizer of
$g$ in $G$ (choose one representative for each class).
Its entries are given by
\begin{equation}\label{quantS}
S_{(\bar a,\chi),(\bar b,\chi')} := \frac{1}{|C_G(a)||C_G(b)|}
\sum_{\substack{g\in G \\ agbg^{-1}=gbg^{-1}a}}
\BAR{\chi(gbg^{-1})}\BAR{\chi'(g^{-1}ag)}, \end{equation}
\[ T_{(\bar a,\chi),(\bar a',\chi')} := \delta_{\bar a,\bar a'}\delta_{\chi,\chi'}
\frac{\chi(a)}{\chi(e)} \]
where $C_G(a)$ is the centralizer of $a$ in $G$ and $e$
is the neutral element of $G$.

Most indecomposable integral modular data known to the author
have fusion algebras isomorphic to the fusion algebra corresponding to the
$S$-matrix of the quantum double of a finite group,
so they have an $S$-matrix which is of the above form
up to an independent permutation of rows and columns.
However, the 2-dimensional modular
data of example \ref{moddatZ2} and the infinite sequence of
proposition \ref{infintexpl} do not have this form.

For example, there are 6 non isomorphic $S$-matrices with
fusion algebras isomorphic to the algebra of the quantum double
of $S_3$. Four of them have $T$-matrices which make them
to modular data, but the groups spanned by $S$ and $T$
fall into two isomorphism classes, the one corresponding
to the quantum double is
$\PSL_2(\Zz/6\Zz)$ and has $72$ elements,
the other one is $\SL_2(\Zz/12\Zz)/N$ where $N$ is a normal subgroup of order $4$
(it has $288$ elements).

The following proposition has been pointed out to the author
by the referee and by Gannon after the present paper had been submitted:

\begin{prop}\label{referee2}
Let $G$ be a finite group and $S$ be the $S$-matrix of the quantum
double of $G$. If $S$ is rational, then $G$ is solvable and the only
possible composition factors are the cyclic groups $Z_2$ and $Z_3$.
\end{prop}
\begin{proof}
This is a consequence of the Galois symmetry of modular data mentioned in \cite{CGR} or \cite{tG}: restricted to the finite group modular data, it says that for any Galois automorphism $\sigma_\ell$ in $\Gal(Q[\zeta_m]/Q)=Z_m^\times$, where $\Qz[\zeta_m]$ is the $m$-th cyclotomic field and $m$ is the exponent of $G$, there is a permutation (also denoted by $\sigma_\ell$) of the indices $1,\ldots,n$ such that
\[ \sigma_\ell(S_{ij})=S_{\sigma_\ell(i),j},\quad
T_{\sigma_\ell(i),\sigma_\ell(i)}=T_{ii}^{\ell^2}.\]
The rationality of $S$ implies that the permutation is trivial on the indices, so $T_{ii}^{\ell^2}=T_{ii}$ for all $\ell$ coprime to $m$ ($T_{ii}$ will be $m$-th roots of unity). This requires that $m$ divides $24$. Now by Burnside's $p^aq^b$ theorem, such a $G$ must be solvable, so the only possible composition factors are $Z_2$ and $Z_3$.
\end{proof}

In the next section, we give some infinite sequences of integral modular data
arising as subdirect products.
The following list contains the groups with at most $100$ elements and rational
$S$-matrix which are not included in
our sequences of subdirect products and are not direct products
of smaller groups. A pair $(a,b)$ corresponds to the $b$-th group
of order $a$ in the database of small groups of
the computer algebra system {\sc GAP}:
\[ (32,49), (32,50), (64,134), (64,137), (64,138), (64,177),
 (64,178), (64,182). \]

\subsection{Subdirect products}

We consider subdirect products with certain properties.
Let $H$ be a finite group, $N$ a normal subgroup and
$\mu : H \rightarrow H/N$ the canonical projection.
The subgroup $G$ of the direct product $H^r$ consisting
of all $(g_1,\ldots,g_r)$ such that $\mu(g_i)=\mu(g_j)$
for all $1\le i,j \le r$ is the subdirect product
$G = H\times_\mu \cdots \times_\mu H$.

\begin{exm}\label{S3Z3}
Take $H=S_3$ and $N=Z_3$.
Then the group $G$ has $\frac{1}{2}(3^r+3)$ conjugacy classes.
The centralizer of an element of order $2$ in $G$ is $Z_2$ and
the centralizer of a non trivial element in $Z_3^r$ is
$Z_3^r$, so the $S$-matrix of the quantum double of $G$ has dimension
$\frac{1}{2}(3^{2r}+7)$.
\end{exm}

We will prove that the $S$-matrix of the quantum double of
subdirect products of the above form is rational in the following
cases:
\begin{enumerate}
\item $H=S_3, D_4, Q_8$ and arbitrary normal subgroup $N$,
\item $H=S_4$ and $N\ne A_4$.
\end{enumerate}

To prove this for the case $H=S_4$, we need a technical
proposition and a computer program (which is still simpler than
determining all conjugacy classes and centralizers of the
subdirect products in details). The other cases may also
be proved using the same technique, but there are also
easier arguments which we do not want to keep back.

The principle is always the same: To each summand in the
formula (\ref{quantS}) for the $S$-matrix, we find another
summand which is its complex conjugate.
This suffices for the
groups whose exponent divides $6$. For $S_4$, we have to
consider more automorphisms.

\subsection{The case $H=D_4$}

For the dihedral group, we use the following lemma:

\begin{lem}\label{subdirnorm}
Let $N\unlhd H$ be finite groups and $G=H\times_\mu\cdots\times_\mu H$
be the subdirect product with respect to $\mu : H \rightarrow H/N$.
Assume that there exists a $\sigma\in H$ such that
\[ \tau^\sigma := \sigma^{-1} \tau \sigma = \tau^{-1} \]
for all $\tau \in N$. Then for the $S$-matrix of the quantum
double of $G$ we have
\[ S_{(\BAR{(a_1,\ldots,a_r)},\chi), (\BAR{(a'_1,\ldots,a'_r)},\chi')} \in \Rz \]
for $(a_1,\ldots,a_r)$, $(a'_1,\ldots,a'_r) \in G$
with $a_i, a'_i \in N$ for all $1\le i \le r$ and
any $\chi\in \Irr(C_G((a_1,\ldots,a_r)))$, $\chi'\in \Irr(C_G((a'_1,\ldots,a'_r)))$.

Further, the entry
$S_{(\BAR{(a_1,\ldots,a_r)},\chi), (\BAR{(a'_1,\ldots,a'_r)},\chi')}$
is rational if $a_1,\ldots,a_r$ lie in $N$,
$a'_1,\ldots,a'_r$ are arbitrary and $\chi'$ is a rational character.
\end{lem}

\begin{proof}
Let $a_i,b_i\in N$ and $g_i\in H$ be such that $a_ig_ib_ig_i^{-1} = g_ib_ig_i^{-1}a_i$
for $1\le i \le r$.
Then for $h_i:=\sigma^{-1} g_i$ we get
\[ h_ib_ih_i^{-1}=g_ib_i^{-1}g_i^{-1}, \quad h_i^{-1}a_ih_i=g_i^{-1}a_i^{-1}g_i,
\quad a_ih_ib_ih_i^{-1}=h_ib_ih_i^{-1}a_i, \]
and if $(g_1,\ldots,g_r)\in G$ then $(h_1,\ldots,h_r)\in G$. So
in the formula for $S$, we can pair the summands
\[ \chi(\underline{g}\underline{b}\underline{g}^{-1})
\chi'(\underline{g}^{-1}\underline{a}\underline{g})
+ \chi(\underline{h}\underline{b}\underline{h}^{-1})
\chi'(\underline{h}^{-1}\underline{a}\underline{h}) \in \Rz \]
where $\underline{g}=(g_1,\ldots,g_r)$ etc.

For the second assertion, use exactly the same argument.
\end{proof}

\begin{rem}\label{ratchar}
If $G$ is a finite group, $|G|=k$ and $\sigma\in G$,
then $\chi(\sigma)$ is rational for every character $\chi$
of $G$ if and only if $\sigma$ is conjugate to $\sigma^m$ for
every integer $m$ prime to $k$.
\end{rem}

\begin{exm}\label{D4Z4}
Let $G=H\times_\mu\cdots\times_\mu H$ be the
subdirect product with $r$ factors
where $H=D_4$ is the dihedral group with $8$ elements
and $\mu$ is the projection $H \rightarrow H/N$ where
$N=Z_4$ is the cyclic normal subgroup of $D_4$ of order $4$.
Let $\eta$ be a generator of $N$.

For all $\sigma\in H\backslash N$ and $\tau\in N$ we have $\sigma^2=1$
and $\tau^\sigma = \tau^{-1}$ (therefore also $\sigma^\tau = \tau^2 \sigma$).
Further, $\sigma^{\sigma'}\in \{\sigma,\eta^2\sigma \}$ for all
$\sigma,\sigma'\in H \backslash N$
and $N\times\ldots\times N\le G$ is abelian and normal in $G$.
So we get $2^r+\frac{1}{2}(4^r-2^r)$ conjugacy classes of $G$ of
the form
\[ \{ (\eta_1,\ldots,\eta_r), (\eta_1,\ldots,\eta_r)^{-1} \} \]
with $\eta_1,\ldots,\eta_r \in N$, and
$2^r$ conjugacy classes of the form
\[ \{ (\eta_1^2 \sigma_1,\ldots,\eta_r^2 \sigma_r)
\mid \eta_1,\ldots,\eta_r \in N \} \]
with $\sigma_1,\ldots,\sigma_r \in H\backslash N$.
Hence $G$ has $2^{r+1}+2^{2r-1}-2^{r-1}$ conjugacy classes.

The centralizer of $(\eta_1,\ldots,\eta_r) \in N\times\ldots\times N$ are:
\[ C_G((\eta_1,\ldots,\eta_r)) = \begin{cases}
G & \mbox{if } \eta_i \in \{1,\eta^2 \} \\
N\times\ldots\times N & \mbox{else.} \end{cases} \]
The centralizer of $(\sigma_1,\ldots,\sigma_r)
\in H\backslash N\times\ldots\times H\backslash N$ are:
\[ C_G((\sigma_1,\ldots,\sigma_r)) = \] \[ = \{ (\tau_1,\ldots,\tau_r) \mid
(\tau_i \in \{1,\eta^2\} \quad \forall i) \mbox{ or }
(\tau_i\in H\backslash N \mbox{ and } \sigma_i^{\tau_i}=\sigma_i \quad \forall i) \} \]
\[ \cong Z_2^{r+1} \]
because $C_G((\sigma_1,\ldots,\sigma_r))$ is a subgroup
of order $2^{r+1}$ of
\[ C_H(\sigma_1)\times\ldots\times C_H(\sigma_r)
\cong (Z_2\times Z_2)^r \]
which is a group of exponent $2$.

The rationality of the $S$-matrix of the quantum
double now follows from lemma \ref{subdirnorm}:

\begin{prop}
The finite group modular data associated to $G$ is integral
for all $r\in\Nz$.
\end{prop}

\begin{proof}
We only need to check that all characters of $G$ are rational,
because then we know that non rational characters
of centralizers only appear in the case of conjugacy classes lying
in $N\times \ldots \times N$. And these are handled by
lemma \ref{subdirnorm}.

But by remark \ref{ratchar}, it suffices to prove that
an element of $G$ is conjugate to all its odd powers
($|G|=2^{2r+1}$).
This follows from $\sigma^3=\sigma$ and $\tau^\sigma=\tau^{-1}$
for $\sigma\in H\backslash N$, $\tau\in N$.
\end{proof}
\end{exm}

\subsection{The case $H=Q_8$, $N\ne H$}
For a finite group $G$, define
\[ k_G(a,b):=|\{ g\in G \mid agbg^{-1} = gbg^{-1}a \}| \]
to be the number of summands in the sum of formula (\ref{quantS})
for the $S$-matrix of the quantum double.

The conjugacy classes are similar to those in Example \ref{D4Z4}.
The number of elements of the centralizer also coincide. There
are two important differences: The elements of $H\backslash N$
do not necessarily have order $2$ and their centralizer are not
necessarily isomorphic to $Z_2$. But the following holds:

\begin{prop}\label{ungl}
Let $H:=Q_8$ and $N\unlhd H$, $N\ne H$ be a normal subgroup.

Then $k_H(a,b)\in\{0,|H|\}$ for all $a,b\in H$ and we have
bijections
\[ \varphi_\sigma : H \rightarrow H \]
for $\sigma\in H$ such that
$\varphi_\sigma(h)\ne h$ for all $h\in H$ and:
\[ \forall b,g \in H, a \in N\sigma \mbox{ with } k_H(a,b)\ne 0 \: : \]
\[ gbg^{-1}= \varphi_\sigma(g)b^{-1}\varphi_\sigma(g)^{-1} \mbox{ and }
g^{-1}ag = \varphi_\sigma(g)^{-1}a^{-1} \varphi_\sigma(g). \]
\end{prop}

\begin{proof}
Compute $k_H(a,b)$ for all $a,b\in H$ and $\varphi_\sigma$ explicitly.
\end{proof}

Now let $G$ be the subdirect product $H\times_\mu\ldots\times_\mu H$
with $r$ factors
where $\mu : H \rightarrow H/N$ is the canonical projection,
$H=Q_8$ and $N$ is a
normal subgroup of $H$, $N \ne H$.
Then by proposition \ref{ungl}, for all $\underline{a},
\underline{b}, \underline{g}\in G$,
$\chi\in\Irr(C_G(\underline{a}))$ and $\chi'\in\Irr(C_G(\underline{b}))$,
\[ \BAR{\chi(\underline{g}\underline{b}\underline{g}^{-1})}
\BAR{\chi'(\underline{g}^{-1}\underline{a}\underline{g})} =
{\chi(\underline{\varphi_\sigma}(\underline{g})\underline{b}
\underline{\varphi_\sigma}(\underline{g})^{-1})}
{\chi'(\underline{\varphi_\sigma}(\underline{g})^{-1} \underline{a}
\underline{\varphi_\sigma}(\underline{g}))} \]
where $\sigma\in H$ is such that $\underline{a}\in N\sigma\times\ldots\times N\sigma$
and $\underline{\varphi_\sigma}:=(\varphi_\sigma,\ldots,\varphi_\sigma)$.
So in the sum of formula (\ref{quantS}), we can always pair $\underline{g}$ with
$\underline{\varphi_\sigma}(\underline{g})$ and obtain a rational number.
We have proved:

\begin{prop}
The finite group modular data associated to $G$ is integral
for all $r\in\Nz$.
\end{prop}

\subsection{The cases $H=S_3$, $N$ arbitrary and $H=S_4$, $N = Z_2\times Z_2$}

For $H=S_3$, the situation is simple; we have almost
done the work in example (\ref{S3Z3}) and therefore leave out a
direct proof here.

More generally,
let $H$ be a finite group and $N\unlhd H$ be a normal subgroup.
For $a,b\in H$, define
\[ K_H(a,b):=\{ g\in H \mid agbg^{-1} = gbg^{-1}a \} \]
to be the set of $g$ appearing in the sum of formula (\ref{quantS})
for the $S$-matrix of the quantum double of $H$.

\begin{prop}\label{phiab}
Assume that for all $m$ coprime to the exponent of $H$,
there exists for all $a,b\in H$ a bijection
\[ \varphi^{(m)}_{a,b} : K_H(a,b) \rightarrow K_H(a,b) \]
such that $\varphi^{(m)}_{a,b}(g)\ne g$ and
\[ gbg^{-1}= \varphi^{(m)}_{a,b}(g)b^m\varphi^{(m)}_{a,b}(g)^{-1}, \quad
g^{-1}ag = \varphi^{(m)}_{a,b}(g)^{-1}a^m \varphi^{(m)}_{a,b}(g) \]
for all $a,b\in H$, $g\in K_H(a,b)$.

Further, assume that for all $a',b'\in H$, $g'\in K_H(a',b')$,
the equalities
\[ Na=Na',\quad Nb=Nb',\quad Ng=Ng' \]
imply $N\varphi^{(m)}_{a,b}(g) = N\varphi^{(m)}_{a',b'}(g')$.

Then the $S$-matrix of the quantum double of the subdirect
product $G:=H\times_\mu \ldots \times_\mu H$, where
$\mu : H \rightarrow H/N$ is the canonical projection,
is rational.
\end{prop}

\begin{proof}
Let $\underline{a}=(a_1,\ldots,a_r),
\underline{b}=(b_1,\ldots,b_r),\underline{g}=(g_1,\ldots,g_r)\in G$
with $g_i\in K_H(a_i,b_i)$, and denote the set of such $\underline{g}$
by $\underline{K}_H(\underline{a},\underline{b})$.
By our assumptions,
\[ \underline{\varphi^{(m)}_{a,b}}:=
(\varphi^{(m)}_{a_1,b_1},\ldots,\varphi^{(m)}_{a_r,b_r}) \]
is a well defined map from $\underline{K}_H(\underline{a},\underline{b})$
to $\underline{K}_H(\underline{a},\underline{b})$ for all $m$ coprime to
the exponent of $H$.
Therefore, for irreducible characters
$\chi\in\Irr(C_G(\underline{a}))$ and $\chi'\in\Irr(C_G(\underline{b}))$,
if the summand $\chi(\underline{g}\underline{b}\underline{g}^{-1})
\chi'(\underline{g}^{-1}\underline{a}\underline{g})$ appears
in formula (\ref{quantS}), then all its Galois conjugate also
appear. Reordering the sum in orbits under Galois conjugation then yields
a sum of rational numbers.

(Notice that since the $\varphi^{(m)}_{a,b}$ are bijections, we may assume
$\varphi^{(m)}_{a,b}=(\varphi^{(m^{-1})}_{a,b})^{-1}$ and hence by the action
of the Galois group we obtain indeed orbits.)
\end{proof}

\begin{prop}\label{S4rat}
For $H=S_3$, $N$ arbitrary or $H=S_4$, $N=Z_2\times Z_2$, the maps
$\varphi^{(m)}_{a,b}$ of Proposition \ref{phiab} exist.
So the $S$-matrix of the quantum double of the
subdirect product is rational.
\end{prop}

\begin{proof}
A solution for the equations given in proposition
\ref{phiab} can be computed by a recursion on a computer.
\end{proof}

\begin{rem}
The maps $\varphi^{(m)}_{a,b}$ of
proposition \ref{phiab} also exist in the cases $H=D_4$,
$H=Q_8$ and $N\unlhd H$ arbitrary.
\end{rem}

\section{Open questions}

The analysis of integral modular data reveals some
difficulties which we will certainly also encounter in
the case of arbitrary modular data. Therefore it
could be interesting to answer the following
questions (some of them could be easy to answer):

\begin{enumerate}
\item Is there an upper bound for the largest norm of
an integral modular datum in a given dimension?
\item What are the integral Fourier matrices with `small'
structure constants (for example in $\{-1,0,1\}$)?
\item Are there integral Fourier matrices with the same
degrees, nonnegative structure constants and non isomorphic fusion algebras
(see example \ref{example2})?
\item Is there an integral Fourier matrix with a unique
norm, negative structure constants and an entry unequal to $\pm 1$?
\item Does an integral Fourier matrix $s$ exist with a column $k$
such that $\forall i$ $s_{ik}\ne 0$ and $\exists i$ $s_{ik} \ne \pm 1$?
\item What are the degrees of congruence multiplication tables?
\item Is there an integral Fourier matrix with only one degree
equal to $1$ and of dimension greater than $1$?
\item Does the integrality of the structure constants of a
congruence multiplication table already follow from the other axioms?
\item What is the `modulus' (or the corresponding ideal) in the
case of arbitrary modular data?
\item Are there modular data with rational $S$ and $T$ matrices which
do not belong to the infinite sequence of proposition \ref{infintexpl} (probably yes)?
\item Classify integral finite group modular data.
\end{enumerate}

\section{Appendix}

\subsection{Computation of all integral Fourier matrices up to dimension 12}

\subsubsection{The norms}

The first step is to list all possible norms. For
a given dimension $n>1$, we have to compute all tuples
$d_1 \ge \ldots \ge d_n$ with the properties
\[ d_j|d_1, \quad \sqrt{\frac{d_1}{d_j}}\in\Zz, \quad
a := \sum_i \frac{1}{d_i} = 1 \]
for all $1\le j\le n$. For this, we start with
$2\le d_n \le n$.

By lemma \ref{oddeven}, if $n$ is odd, we may assume that all
$d_i$ are squares. If $n$ is even, then the square
free part of all $d_i$ is the same, so in this case
we divide $d_n$ by its square free part $p$ and start
again in dimension $n-1$ and the above sum now
sums up to $a:=p-\frac{p}{d_n}$; we replace $d_n$ by
$\frac{d_n}{p}$ which is a square.

Then we start a recursion. There are two lower bounds
for the next $d_i$ in each step given by the above
number $a$ (from which we subtract the inverse of
the last $d_i$ in each step)
and by $d_{i+1}$. The upper bound is $\frac{n}{a}$.

Solutions $d_1,\ldots,d_n$ to the equation
$\frac{1}{d_1}+\ldots+\frac{1}{d_n}=1$ are called Egyptian fractions.
For a given $n$, the $d_i$ are bounded by the $n$-th term
of Sylvester's sequence: $a_{n+1}=a_n^2-a_n+1$, $a_1=2$ (\cite{IS}
gives the reference \cite{dS}).
This information is not useful for the computation because
the numbers $a_n$ are very large.

\subsubsection{Finding the matrices to a given norm vector}

The second step is much more complicated, because there
is a lot of information to exploit in the algorithm.
Again, we construct the matrix in a recursive way,
row by row.
Because it is symmetric, we still have to compute
a square matrix after each step.

There are two important observations: it is very useful
to test if the structure constants are integers after
each step (only for those structure constants which may
already be computed of course).
Also, if there is a norm with high multiplicity, say $m$, then
the recursion will compute almost all permutations of the correcponding
block, so approximatively $m!$ matrices. So it is
crucial to keep only representatives of the matrices
up to simultaneous permutations of rows and columns before
the next recursion step.
It is also possible to include several bounds into the computation.

Finding all matrices up to dimension $11$ takes just a few minutes,
whereas dimension $12$ takes about a week on an ordinary PC.

Finding congruence multiplication tables is of course much
easier, but you still need to compute all norm vectors for one dimension.
This begins to become difficult at dimension $13$, even when using
the information of the above propositions.

\subsubsection{Finding $T$-matrices}

For a given $s$-matrix $s$, a $T$-matrix
\[ T=\diag(t_1,\ldots,t_n) \]
is a solution of the equations
\begin{equation}\label{tmsol}
\sum_k t_i t_j t_k \frac{s_{ik} s_{kj}}{\sqrt{d_k}} = s_{ij}
\end{equation}
for all $1\le i,j\le n$.
If $T$ has finite order, then $t_1,\ldots,t_n$ are roots of unity.
In particular, for $i=j=1$ we get
$\sum_k \frac{\sqrt{d_1}}{d_k} t_1^2 t_k =1$.
And if $s_{ij} = 0 $ for some $1\le i,j\le n$, then
$\sum_k s_{ik} s_{jk} t_k = 0$.

If $s$ has small dimension, then there is no need
to elaborate a complicated algorithm: it suffices to compute
the Gr\"obner basis of the system of equations (\ref{tmsol})
in an adequate field. All matrices of the appendix are small
enough for this naive method.

\subsection{An integral modular datum in dimension 16}

The following matrix is the $S$-matrix of an integral
modular datum. A corresponding $T$-matrix does not have
finite order, so it is not a modular datum in the sense
of Gannon. But the structure constants are nonnegative.
Remark also, that it is a congruence multiplication table.

\begin{tiny}
\[ \left( \begin{array}{cccccccccccccccc}
1 & 1 & 2 & 6 & 6 & 6 & 6 & 6 & 6 & 6 & 6 & 14 & 14 & 14 & 21 & 21 \\
1 & 1 & 2 & 6 & 6 & 6 & 6 & 6 & 6 & 6 & 6 & 14 & 14 & 14 & -21 & -21 \\
1 & 1 & 2 & 6 & 6 & 6 & 6 & 6 & 6 & 6 & 6 & -7 & -7 & -7 & 0 & 0 \\
1 & 1 & 2 & 6 & -1 & -1 & -1 & -1 & -1 & -1 & -1 & 0 & 0 & 0 & 0 & 0 \\
1 & 1 & 2 & -1 & 6 & -1 & -1 & -1 & -1 & -1 & -1 & 0 & 0 & 0 & 0 & 0 \\
1 & 1 & 2 & -1 & -1 & 6 & -1 & -1 & -1 & -1 & -1 & 0 & 0 & 0 & 0 & 0 \\
1 & 1 & 2 & -1 & -1 & -1 & 6 & -1 & -1 & -1 & -1 & 0 & 0 & 0 & 0 & 0 \\
1 & 1 & 2 & -1 & -1 & -1 & -1 & 6 & -1 & -1 & -1 & 0 & 0 & 0 & 0 & 0 \\
1 & 1 & 2 & -1 & -1 & -1 & -1 & -1 & 6 & -1 & -1 & 0 & 0 & 0 & 0 & 0 \\
1 & 1 & 2 & -1 & -1 & -1 & -1 & -1 & -1 & -1 & 6 & 0 & 0 & 0 & 0 & 0 \\
1 & 1 & 2 & -1 & -1 & -1 & -1 & -1 & -1 & 6 & -1 & 0 & 0 & 0 & 0 & 0 \\
1 & 1 & -1 & 0 & 0 & 0 & 0 & 0 & 0 & 0 & 0 & 2 & -1 & -1 & 0 & 0 \\
1 & 1 & -1 & 0 & 0 & 0 & 0 & 0 & 0 & 0 & 0 & -1 & -1 & 2 & 0 & 0 \\
1 & 1 & -1 & 0 & 0 & 0 & 0 & 0 & 0 & 0 & 0 & -1 & 2 & -1 & 0 & 0 \\
1 & -1 & 0 & 0 & 0 & 0 & 0 & 0 & 0 & 0 & 0 & 0 & 0 & 0 & -1 & 1 \\
1 & -1 & 0 & 0 & 0 & 0 & 0 & 0 & 0 & 0 & 0 & 0 & 0 & 0 & 1 & -1 \\
\end{array} \right) \]
\end{tiny}

A solution for a $T$-matrix is for example
\[ T:=\diag(1,1,1,1,1,1,1,1,1,\alpha,\alpha^{-1},1,\zeta,\zeta^2,1,-1),\]
where $\alpha$ is a root of $X^2+5X+1$ and $\zeta$ a third root
of unity.

\subsection{A `rational' modular datum in dimension $18$}

The following matrix is the $S$-matrix of an integral
modular datum. All corresponding $T$-matrices have
finite order and the structure constants are nonnegative,
so it is a modular datum in the sense of Gannon.
It is not a congruence multiplication table.

\begin{tiny}
\[ \left( \begin{array}{cccccccccccccccccc}
1 & 1 & 1 & 1 & 3 & 3 & 3 & 3 & 3 & 3 & 3 & 3 & 3 & 3 & 3 & 3 & 4 & 4 \\
1 & 1 & 1 & 1 & 3 & 3 & 3 & 3 & -3 & -3 & -3 & -3 & -3 & -3 & -3 & -3 & 4 & 4 \\
1 & 1 & 1 & 1 & -3 & -3 & -3 & -3 & 3 & 3 & 3 & 3 & -3 & -3 & -3 & -3 & 4 & 4 \\
1 & 1 & 1 & 1 & -3 & -3 & -3 & -3 & -3 & -3 & -3 & -3 & 3 & 3 & 3 & 3 & 4 & 4 \\
1 & 1 & -1 & -1 & -1 & 1 & 1 & -1 & 1 & 1 & -1 & -1 & 1 & 1 & -1 & -1 & 0 & 0 \\
1 & 1 & -1 & -1 & 1 & -1 & -1 & 1 & 1 & 1 & -1 & -1 & -1 & -1 & 1 & 1 & 0 & 0 \\
1 & 1 & -1 & -1 & 1 & -1 & -1 & 1 & -1 & -1 & 1 & 1 & 1 & 1 & -1 & -1 & 0 & 0 \\
1 & 1 & -1 & -1 & -1 & 1 & 1 & -1 & -1 & -1 & 1 & 1 & -1 & -1 & 1 & 1 & 0 & 0 \\
1 & -1 & 1 & -1 & 1 & 1 & -1 & -1 & -1 & 1 & 1 & -1 & 1 & -1 & 1 & -1 & 0 & 0 \\
1 & -1 & 1 & -1 & 1 & 1 & -1 & -1 & 1 & -1 & -1 & 1 & -1 & 1 & -1 & 1 & 0 & 0 \\
1 & -1 & 1 & -1 & -1 & -1 & 1 & 1 & 1 & -1 & -1 & 1 & 1 & -1 & 1 & -1 & 0 & 0 \\
1 & -1 & 1 & -1 & -1 & -1 & 1 & 1 & -1 & 1 & 1 & -1 & -1 & 1 & -1 & 1 & 0 & 0 \\
1 & -1 & -1 & 1 & 1 & -1 & 1 & -1 & 1 & -1 & 1 & -1 & -1 & 1 & 1 & -1 & 0 & 0 \\
1 & -1 & -1 & 1 & 1 & -1 & 1 & -1 & -1 & 1 & -1 & 1 & 1 & -1 & -1 & 1 & 0 & 0 \\
1 & -1 & -1 & 1 & -1 & 1 & -1 & 1 & 1 & -1 & 1 & -1 & 1 & -1 & -1 & 1 & 0 & 0 \\
1 & -1 & -1 & 1 & -1 & 1 & -1 & 1 & -1 & 1 & -1 & 1 & -1 & 1 & 1 & -1 & 0 & 0 \\
1 & 1 & 1 & 1 & 0 & 0 & 0 & 0 & 0 & 0 & 0 & 0 & 0 & 0 & 0 & 0 & 1 & -2 \\
1 & 1 & 1 & 1 & 0 & 0 & 0 & 0 & 0 & 0 & 0 & 0 & 0 & 0 & 0 & 0 & -2 & 1 \\
\end{array} \right) \]
\end{tiny}

A solution for a $T$-matrix is for example
\[ T:=\diag(-1, -1, -1, -1, -1, 1, 1, -1, -1, 1, 1, -1, 1, -1, -1, 1, \zeta^{-1}, \zeta),\]
where $\zeta$ a sixth root of unity. The group spanned
by $S$ and $T$ is $\PSL_2(\Zz/6\Zz)$. Permuting columns of $S$ yields
another Fourier matrix which has $T$-matrices such that the groups
spanned by $S$ and $T$ are isomorphic to $\SL_2(\Zz/12\Zz)/N$ where
$N$ is a normal subgroup of order $4$. Remark that these groups
are those that appear for the fusion algebra corresponding
to the quantum double of $S_3$.

As mentioned in the introduction, this modular data is
similar to one of the modular data corresponding to
the exceptional complex reflection group $G_8$. Section $8.2.3$
of \cite{mC} describes a construction for its fusion algebra
as a factor ring of the quantum double of the $9$-th group of
order $36$ in the numbering of {\sc GAP}.

\subsection{Dimension 1 to 12}

This is a list of all isomorphism classes of matrices $s$
satisfying the axioms (1), (2), (3), (4) of definition \ref{intmoddef}
up to dimension $12$.

A permutation $\sigma$ as an index means that it is
a non isomorphic matrix obtained by
permuting the columns of this matrix by $\sigma$.
If there are two permutations, then the first one is for the
columns, the second one for the rows.

The exponent $+6$ stands for the construction of
proposition \ref{oldcmtp6}.

\begin{tiny}
\[ s(2,1) := \left( \begin{array}{cc}
1 & 1 \\
1 & -1 \\
\end{array} \right), \quad
s(4,1) := \left( \begin{array}{cccc}
1 & 1 & 1 & 1 \\
1 & 1 & -1 & -1 \\
1 & -1 & 1 & -1 \\
1 & -1 & -1 & 1 \\
\end{array} \right), \quad
s(4,1)_{(3, 4)}, \quad
s(7,1) := \left( \begin{array}{c}
1 \\
\end{array} \right)^{+6} \]
\[ s(8,1) :=
s(2,1)^{\otimes 3}, \quad
s(8,2) := \left( \begin{array}{cccccccc}
1 & 1 & 2 & 2 & 2 & 2 & 3 & 3 \\
1 & 1 & 2 & 2 & 2 & 2 & -3 & -3 \\
1 & 1 & -1 & -1 & -1 & 2 & 0 & 0 \\
1 & 1 & -1 & -1 & 2 & -1 & 0 & 0 \\
1 & 1 & -1 & 2 & -1 & -1 & 0 & 0 \\
1 & 1 & 2 & -1 & -1 & -1 & 0 & 0 \\
1 & -1 & 0 & 0 & 0 & 0 & 1 & -1 \\
1 & -1 & 0 & 0 & 0 & 0 & -1 & 1 \\
\end{array} \right) \]
\[ s(8,2)_{(7, 8)}, \quad
 s(8,2)_{(4, 5)}, \quad
 s(8,2)_{(4, 5)(7, 8)}, \quad
 s(8,2)_{(3, 6)(4, 5)}, \quad
 s(8,2)_{(3, 6)(4, 5)(7, 8)} \]
\[ s(9,1) := \left( \begin{array}{ccc}
1 & 2 & 2 \\
1 & \frac{1}{2} & -1 \\
1 & -1 & \frac{1}{2} \\
\end{array} \right)^{+6}, \quad
s(9,1)_{(5, 6)}, \quad
s(9,2) := \left( \begin{array}{ccc}
1 & \frac{1}{2} & 1 \\
1 & 2 & -2 \\
1 & -1 & -\frac{1}{2} \\
\end{array} \right)^{+6} \]
\[ s(10,1) := \left( \begin{array}{cccccccccc}
1 & 1 & 1 & 1 & 1 & 2 & 3 & 3 & 3 & 6 \\
1 & 1 & 1 & 1 & 1 & 2 & 3 & 3 & 3 & -6 \\
1 & 1 & 1 & 1 & -5 & -4 & 3 & 3 & -3 & 0 \\
1 & 1 & 1 & -5 & 1 & -4 & 3 & -3 & 3 & 0 \\
1 & 1 & -5 & 1 & 1 & -4 & -3 & 3 & 3 & 0 \\
1 & 1 & -2 & -2 & -2 & 2 & 0 & 0 & 0 & 0 \\
1 & 1 & 1 & 1 & -1 & 0 & -1 & -1 & 1 & 0 \\
1 & 1 & 1 & -1 & 1 & 0 & -1 & 1 & -1 & 0 \\
1 & 1 & -1 & 1 & 1 & 0 & 1 & -1 & -1 & 0 \\
1 & -1 & 0 & 0 & 0 & 0 & 0 & 0 & 0 & 0 \\
\end{array} \right), \quad
s(10,1)_{(7, 9), (3, 5)} \] 
\[ s(10,2) := s(4,1)^{+6}, \quad
s(10,2)_{(6, 7)} \] 
\[ s(10,3) := \left( \begin{array}{cccccccccc}
1 & 1 & 2 & 3 & 3 & 4 & 4 & 4 & 6 & 6 \\
1 & 1 & 2 & 3 & 3 & 4 & 4 & 4 & -6 & -6 \\
1 & 1 & 2 & 3 & 3 & -2 & -2 & -2 & 0 & 0 \\
1 & 1 & 2 & -1 & -1 & 0 & 0 & 0 & 2 & -2 \\
1 & 1 & 2 & -1 & -1 & 0 & 0 & 0 & -2 & 2 \\
1 & 1 & -1 & 0 & 0 & 1 & 1 & -2 & 0 & 0 \\
1 & 1 & -1 & 0 & 0 & 1 & -2 & 1 & 0 & 0 \\
1 & 1 & -1 & 0 & 0 & -2 & 1 & 1 & 0 & 0 \\
1 & -1 & 0 & 1 & -1 & 0 & 0 & 0 & 0 & 0 \\
1 & -1 & 0 & -1 & 1 & 0 & 0 & 0 & 0 & 0 \\
\end{array} \right), \quad
s(10,3)_{(6, 8)} \] 
\[ s(11,1) := \left( \begin{array}{ccccc}
1 & 2 & 2 & 6 & 6 \\
1 & 2 & 2 & \frac{3}{2} & -3 \\
1 & 2 & 2 & -3 & \frac{3}{2} \\
1 & \frac{1}{2} & -1 & 0 & 0 \\
1 & -1 & \frac{1}{2} & 0 & 0 \\
\end{array} \right)^{+6}, \quad
s(11,2) := \left( \begin{array}{ccccc}
1 & 1 & 3 & 3 & 4 \\
1 & 1 & -3 & -3 & 4 \\
1 & -1 & 1 & -1 & 0 \\
1 & -1 & -1 & 1 & 0 \\
1 & 1 & 0 & 0 & -\frac{1}{2} \\
\end{array} \right)^{+6} \]
\[ s(11,2)_{(6, 7)}, \quad
s(12,1) := \left( \begin{array}{cccccc}
1 & 1 & 4 & 9 & 9 & 12 \\
1 & 1 & 4 & -9 & -9 & 12 \\
1 & 1 & 4 & 0 & 0 & -\frac{3}{2} \\
1 & -1 & 0 & 1 & -1 & 0 \\
1 & -1 & 0 & -1 & 1 & 0 \\
1 & 1 & -\frac{1}{2} & 0 & 0 & 0 \\
\end{array} \right)^{+6}, \quad
s(12,1)_{(7, 8)} \] 
\[ s(12,2) := \left( \begin{array}{cccccccccccc}
1 & 1 & 2 & 2 & 2 & 2 & 6 & 6 & 6 & 6 & 9 & 9 \\
1 & 1 & 2 & 2 & 2 & 2 & 6 & 6 & 6 & 6 & -9 & -9 \\
1 & 1 & 2 & 2 & 2 & 2 & -3 & -3 & -3 & 6 & 0 & 0 \\
1 & 1 & 2 & 2 & 2 & 2 & -3 & -3 & 6 & -3 & 0 & 0 \\
1 & 1 & 2 & 2 & 2 & 2 & -3 & 6 & -3 & -3 & 0 & 0 \\
1 & 1 & 2 & 2 & 2 & 2 & 6 & -3 & -3 & -3 & 0 & 0 \\
1 & 1 & -1 & -1 & -1 & 2 & 0 & 0 & 0 & 0 & 0 & 0 \\
1 & 1 & -1 & -1 & 2 & -1 & 0 & 0 & 0 & 0 & 0 & 0 \\
1 & 1 & -1 & 2 & -1 & -1 & 0 & 0 & 0 & 0 & 0 & 0 \\
1 & 1 & 2 & -1 & -1 & -1 & 0 & 0 & 0 & 0 & 0 & 0 \\
1 & -1 & 0 & 0 & 0 & 0 & 0 & 0 & 0 & 0 & 1 & -1 \\
1 & -1 & 0 & 0 & 0 & 0 & 0 & 0 & 0 & 0 & -1 & 1 \\
\end{array} \right), \quad
s(12,2)_{(11, 12)} \] 
\[ s(12,3) := \left( \begin{array}{cccccccccccc}
1 & 1 & 2 & 3 & 3 & 6 & 6 & 8 & 8 & 8 & 12 & 12 \\
1 & 1 & 2 & 3 & 3 & 6 & 6 & 8 & 8 & 8 & -12 & -12 \\
1 & 1 & 2 & 3 & 3 & 6 & 6 & -4 & -4 & -4 & 0 & 0 \\
1 & 1 & 2 & 3 & 3 & -2 & -2 & 0 & 0 & 0 & 4 & -4 \\
1 & 1 & 2 & 3 & 3 & -2 & -2 & 0 & 0 & 0 & -4 & 4 \\
1 & 1 & 2 & -1 & -1 & 2 & -2 & 0 & 0 & 0 & 0 & 0 \\
1 & 1 & 2 & -1 & -1 & -2 & 2 & 0 & 0 & 0 & 0 & 0 \\
1 & 1 & -1 & 0 & 0 & 0 & 0 & -1 & -1 & 2 & 0 & 0 \\
1 & 1 & -1 & 0 & 0 & 0 & 0 & -1 & 2 & -1 & 0 & 0 \\
1 & 1 & -1 & 0 & 0 & 0 & 0 & 2 & -1 & -1 & 0 & 0 \\
1 & -1 & 0 & 1 & -1 & 0 & 0 & 0 & 0 & 0 & 0 & 0 \\
1 & -1 & 0 & -1 & 1 & 0 & 0 & 0 & 0 & 0 & 0 & 0 \\
\end{array} \right) \]
\[ s(12,3)_{(8, 10)}, \quad
s(12,3)_{(6, 7)}, \quad
s(12,3)_{(6, 7)(8, 10)} \] 
\[ s(12,4) := \left( \begin{array}{cccccccccccc}
1 & 1 & 2 & 6 & 6 & 6 & 6 & 10 & 10 & 10 & 15 & 15 \\
1 & 1 & 2 & 6 & 6 & 6 & 6 & 10 & 10 & 10 & -15 & -15 \\
1 & 1 & 2 & 6 & 6 & 6 & 6 & -5 & -5 & -5 & 0 & 0 \\
1 & 1 & 2 & 1 & 1 & 1 & -4 & 0 & 0 & 0 & 0 & 0 \\
1 & 1 & 2 & 1 & 1 & -4 & 1 & 0 & 0 & 0 & 0 & 0 \\
1 & 1 & 2 & 1 & -4 & 1 & 1 & 0 & 0 & 0 & 0 & 0 \\
1 & 1 & 2 & -4 & 1 & 1 & 1 & 0 & 0 & 0 & 0 & 0 \\
1 & 1 & -1 & 0 & 0 & 0 & 0 & 1 & 1 & -2 & 0 & 0 \\
1 & 1 & -1 & 0 & 0 & 0 & 0 & 1 & -2 & 1 & 0 & 0 \\
1 & 1 & -1 & 0 & 0 & 0 & 0 & -2 & 1 & 1 & 0 & 0 \\
1 & -1 & 0 & 0 & 0 & 0 & 0 & 0 & 0 & 0 & 1 & -1 \\
1 & -1 & 0 & 0 & 0 & 0 & 0 & 0 & 0 & 0 & -1 & 1 \\
\end{array} \right) \]
\[ s(12,4)_{(11, 12)}, \quad 
 s(12,4)_{(8, 10)}, \quad
 s(12,4)_{(8, 10)(11, 12)}, \quad
 s(12,4)_{(5, 6)} \]
\[ s(12,4)_{(5, 6)(11, 12)} , \quad
 s(12,4)_{(5, 6)(8, 10)} , \quad
 s(12,4)_{(5, 6)(8, 10)(11, 12)} , \quad
 s(12,4)_{(4, 7)(5, 6)} \]
\[ s(12,4)_{(4, 7)(5, 6)(11, 12)} , \quad
 s(12,4)_{(4, 7)(5, 6)(8, 10)} , \quad
 s(12,4)_{(4, 7)(5, 6)(8, 10)(11, 12)} \] 
\[ s(12,5) := \left( \begin{array}{cccccc}
1 & \frac{1}{2} & \frac{1}{2} & \frac{1}{2} & \frac{1}{2} & \frac{1}{2} \\
1 & -1 & -1 & -1 & -1 & 2 \\
1 & -1 & -1 & -1 & 2 & -1 \\
1 & -1 & -1 & 2 & -1 & -1 \\
1 & -1 & 2 & -1 & -1 & -1 \\
1 & 2 & -1 & -1 & -1 & -1 \\
\end{array} \right)^{+6}, \quad
s(12,5)_{(6, 8)} , \quad
 s(12,5)_{(5, 9)(6, 8)} \] 
\[ s(12,6) := \left( \begin{array}{cccccc}
1 & 1 & 1 & 1 & 1 & 2 \\
1 & 1 & 1 & 1 & -2 & -1 \\
1 & 1 & 1 & -2 & 1 & -1 \\
1 & 1 & -2 & 1 & 1 & -1 \\
1 & -2 & 1 & 1 & 1 & -1 \\
1 & -\frac{1}{2} & -\frac{1}{2} & -\frac{1}{2} & -\frac{1}{2} & \frac{1}{2} \\
\end{array} \right)^{+6}, \quad
s(12,6)_{(6, 7)} , \quad
 s(12,6)_{(5, 8)(6, 7)} \] 
\[ s(12,7) := \left( \begin{array}{cccccc}
1 & \frac{1}{2} & 1 & \frac{3}{2} & \frac{3}{2} & \frac{3}{2} \\
1 & \frac{7}{2} & 4 & -\frac{3}{2} & -\frac{3}{2} & -\frac{3}{2} \\
1 & 2 & -2 & 0 & 0 & 0 \\
1 & -\frac{1}{2} & 0 & \frac{1}{2} & \frac{1}{2} & -\frac{3}{2} \\
1 & -\frac{1}{2} & 0 & \frac{1}{2} & -\frac{3}{2} & \frac{1}{2} \\
1 & -\frac{1}{2} & 0 & -\frac{3}{2} & \frac{1}{2} & \frac{1}{2} \\
\end{array} \right)^{+6}, \quad
s(12,7)_{(7, 9)} \] 
\end{tiny}

\subsection*{Acknowledgements}

The author thanks G. Malle for his support and many helpful suggestions.
He is also very grateful to T. Gannon and the referee who have contributed
Proposition \ref{referee2} and to the referee for Proposition \ref{referee1}.

\bibliographystyle{amsplain}
\bibliography{references}

\end{document}